%% file: variety-March-30.tex
\documentclass[11pt]{article}

\usepackage{graphicx,tikz,color,url}
\usepackage{amsmath,amssymb,latexsym}
\usepackage{mathtools}
\usepackage{xspace}
\usepackage[T1]{fontenc}

\newtheorem{theorem}{Theorem}[section]

\newtheorem{proposition}[theorem]{Proposition}

\newtheorem{corollary}[theorem]{Corollary}

\newtheorem{property}[theorem]{Property}

\newcommand{\proof}{\noindent{\bf Proof.\ }}
\newcommand{\qed}{\hfill $\square$ \bigskip}
\newcommand{\cp}{\,\square\,}

\newcommand{\diam}{{\rm diam}}

\newcommand{\mut}{\mu_{\rm t}}
\newcommand{\mud}{\mu_{\rm d}}
\newcommand{\muo}{\mu_{\rm o}}

\definecolor{auburn}{rgb}{0.43, 0.21, 0.1}
\definecolor{blush}{rgb}{0.87, 0.36, 0.51}

\begin{document}

\title{Variety of mutual-visibility problems in graphs}

\author{
Serafino Cicerone $^{a}$ \thanks{Email: \texttt{serafino.cicerone@univaq.it}}
\and
Gabriele Di Stefano $^{a} $\thanks{Email: \texttt{gabriele.distefano@univaq.it}}
\and
Lara Dro\v{z}\dj ek $^{b} $\thanks{Email: \texttt{lara.drozdek@um.si}}
\and
Jaka Hed\v{z}et $^{b,c} $\thanks{Email: \texttt{jaka.hedzet@imfm.si}}
\and
Sandi Klav\v{z}ar $^{d,c,b} $\thanks{Email: \texttt{sandi.klavzar@fmf.uni-lj.si}}
\and
Ismael G. Yero $^{e} $\thanks{Email: \texttt{ismael.gonzalez@uca.es}}
}

\maketitle
\begin{center}
\small
	$^a$ Department of Information Engineering, Computer Science, and Mathematics,
	     University of L'Aquila, Italy \\
	\medskip
    $^b$ Faculty of Natural Sciences and Mathematics, University of Maribor, Slovenia\\
    \medskip
	$^c$ Institute of Mathematics, Physics and Mechanics, Ljubljana, Slovenia\\
	\medskip
	$^d$ Faculty of Mathematics and Physics, University of Ljubljana, Slovenia\\
	\medskip		
	$^e$ Departamento de Matem\'aticas, Universidad de C\'adiz, Spain \\
\end{center}

\begin{abstract}
If $X$ is a subset of vertices of a graph $G$, then vertices $u$ and $v$ are $X$-visible if there exists a shortest $u,v$-path $P$ such that $V(P)\cap X \subseteq \{u,v\}$. If each two vertices from $X$  are $X$-visible, then $X$ is a mutual-visibility set. The mutual-visibility number of $G$ is the cardinality of a largest mutual-visibility set of $G$ and has been already investigated. In this paper a variety of mutual-visibility problems is introduced based on which natural pairs of vertices are required to be $X$-visible. This yields the total, the dual, and the outer mutual-visibility numbers. 
We first show that these graph invariants are related to each other and to the classical mutual-visibility number, and then we prove that the three newly introduced mutual-visibility problems are computationally difficult. According to this result, we compute or bound their values for several graphs classes that include for instance grid graphs and tori. We conclude the study by presenting some inter-comparison between the values of such parameters, which is based on the computations we made for some specific families. 
\end{abstract}


\noindent
AMS Subj.\ Class.\ (2020):  05C12, 05C69, 05C76, 68Q25

\maketitle

\section{Introduction}

Vertex visibility in graphs with respect to a set of vertices has been recently introduced and studied in the sense of the existence of a shortest path between two vertices not containing a third vertex from such set. The visibility property is then understood as a kind of non existence of ``obstacles'' between the two vertices in the mentioned shortest path, which makes them ``visible'' to each other.

Visibility problems in networks have recently attracted the attention of several investigations dealing, in one hand, with theoretical problems arising in the area of graph theory and combinatorics, and in a second hand, with practical problems appearing in the area of computer science. Concerning this latter research line, several contributions focused on applications of visibility problems in some robot navigation models were presented for instance in the works~\cite{aljohani-2018a,bhagat-2020,Cicerone-2023+,diluna-2017,poudel-2021}. On the other hand, in connection with the theoretical studies, the article~\cite{DiStefano-2022} was an introductory contribution that was further continued in~\cite{Cicerone-2023, Cicerone-2022+, tian-2023+}.

In the four last mentioned works, among other results, several contributions were presented aimed to find or bound the largest cardinality of a set of vertices in a graph satisfying that their vertices are pairwise ``visible'', under the assumption of visibility previously commented. The contributions presented in these three works showed several interesting connections like that one existing between such visibility problem and one instance of the very well known Zarankiewicz problem (see \cite[Corollary 3.7]{Cicerone-2023}). The mutual-visibility is also related to the general position problem in graphs;~\cite{klavzar-2023, klavzar-2022, klavzar-2021, patkos-2020, tian-2021, tuite-2022, yao-2022} is a selection of related recent papers, see also references therein.

Results from previous papers on visibility problems in graphs required the use of several powerful tools which introduced modifications on the visibility properties to be taken into account. For instance, the article~\cite{Cicerone-2023} considered vertex visibility in which the vertices in question were also satisfying an independency (the non existence of edges) property between them. This turned out to be very useful whilst considering the visibility problem in Cartesian product graphs. In a similar sense, in~\cite{Cicerone-2022+}, the study of visibility sets in the strong product of graphs led to a modification in the visibility property that was requiring such property to be satisfied not only between the vertices of the set but between every two vertices of the graph.

Research to date showed, among other things, the richness of the topic while considering different styles of visibility situations. Consequently, it seems natural to consider studying the variety of the possible visibility situations that can appear in the investigation. The main purpose of this article is to introduce and motivate the variety. For a given graph $G$, the variety consists of
the mutual-visibility number $\mu(G)$,
the total mutual-visibility number $\mut(G)$,
the outer mutual-visibility number $\muo(G)$, and
the dual mutual-visibility number $\mud(G)$.

In the next section, we first list definitions needed, then formally introduce the variety, and finally provide some basic properties. Section \ref{sec:complexity} is dedicated to computational issues, and in particular to prove that computing each of the three new parameters of the variety is an NP-hard problem. 
In Section~\ref{sec:grids} we determine the value of these invariants for grid-like graphs, that is, Cartesian products of paths or of cycles. In Section~\ref{sec:compare} the invariants are compared and as a byproducts additional exact values determined. We end our exposition with a concluding section, in which we present some possible future research lines that might be of interest for the research community.

\section{Preliminaries and the variety}
\label{sec:preliminaries}

We consider undirected graphs and unless otherwise stated, all graphs in the paper are connected. Given a graph $G$, $V(G)$ and $E(G)$ are used to denote its vertex set and its edge set, respectively. The order of $G$, that is $|V(G)|$, is denoted by $n(G)$. By $m(G)$ we denote the number of edges of $G$, that is $|E(G)|$. If $X\subseteq V(G)$, then $G[X]$ denotes the subgraph of $G$ induced by $X$. The minimum degree of $G$ is denoted by $\delta(G)$. 

The \emph{complement} of a graph $G$ is the graph $\overline{G}$ on the same vertices such that two distinct vertices of $\overline{G}$ are adjacent if and only if they are not adjacent in $G$. The distance function $d_G$ on a graph $G$ is the usual shortest-path distance. The {\em diameter} $\diam(G)$ of $G$ is the maximum distance between pairs of vertices of the graph. A subgraph $G'$ of a graph $G$ is  \emph{convex}, if for every two vertices of $G'$, every shortest path in $G$ between them lies completely in $G'$. A \emph{universal vertex} is a vertex that is adjacent to all other vertices of the graph.  

The \emph{complete graph} (or \emph{clique}) $K_n$, $n\ge 1$, is the graph with $n$ vertices where each pair of distinct vertices are adjacent. For a natural number $n$, we set $[n] = \{1,\ldots, n\}$. The \emph{path graph} $P_n$, $n\ge 2$, is the graph with $V(P_n) = [n]$ such that $i$ is adjacent to $j$ if and only if $|i-j|=1$. The \emph{grid graph} $P_{n} \cp P_{m}$  is the Cartesian product of the paths $P_n$ and $P_m$, that is, $V(P_{n}\cp P_{m}) = \{(i,j):\ i\in [n], j\in [m]\}$ and $(i,j)(k,\ell)\in E(P_{n}\cp P_{m})$ whenever $|i-k|+|j-\ell|=1$.
The \emph{cycle graph} $C_n$, $n\ge 3$, is the graph with $V(C_n) =  [n]$ such that $i$ is adjacent to $j$ if and only if $|i-j|=1$ or $|i-j|=n-1$. The \emph{torus graph} $C_{n} \cp C_{m}$ is the Cartesian product of the cycles $C_n$ and $C_m$, that is, $V(C_{n}\cp C_{m}) = \{(i,j):\ i\in [n], j\in [m]\}$ and $(i,j)(k,\ell)\in E(C_{n}\cp C_{m})$ whenever one of the following conditions holds: $(a)$ $|i-k|+|j-\ell|=1$, or $(b)$ $i=k$ and $|j-\ell|=n-1$, or $(c)$ $|i-k|=n-1$ and $j=\ell$.
A {\em layer} in $P_n\cp P_m$ or $C_n\cp C_m$ is a subgraph induced by the vertices in which one of the coordinates is fixed. Note that each layer is isomorphic either to $P_n$ or $P_m$ in $P_n\cp P_m$, and to $C_n$ or $C_m$ in $C_n\cp C_m$. For $i\in [n]$ and for $j\in [m]$, the corresponding layers will be denoted by $P_m^{(i)}$ and $P_n^{(j)}$ in the grid graph, and by $C_m^{(i)}$ and $C_n^{(j)}$ in the torus graph.

\subsection{The variety}
\label{sec:variety}

The concept of mutual-visibility in graphs was introduced to the literature in~\cite{DiStefano-2022}. While investigating it in the strong product of graphs, the concept of total mutual-visibility has proved to be a natural and necessary tool for its exploration~\cite{Cicerone-2022+}. This has encouraged us to introduce a natural variety of the mutual-visibility in graphs as follows.

Let $G = (V(G), E(G))$ be a graph and $X\subseteq V(G)$. Vertices $u,v\in V(G)$ are {\em $X$-visible} if there exists a shortest $u,v$-path (also called {\em geodesic}) $P$ such that $V(P)\cap X \subseteq \{u,v\}$. Note that each pair of adjacent vertices is $X$-visible. Set $\overline{X} = V(G)\setminus X$. Then we say that $X$ is a
\begin{itemize}
\item \emph{mutual-visibility set}, if every $u,v\in X$ are $X$-visible,
\item \emph{total mutual-visibility set}, if every $u,v\in V(G)$ are $X$-visible,
\item \emph{outer mutual-visibility set}, if every $u,v\in X$ are $X$-visible, and every $u\in X$, $v\in \overline{X}$ are $X$-visible,
\item \emph{dual mutual-visibility set}, if every $u,v\in X$ are $X$-visible, and every $u,v\in \overline{X}$ are $X$-visible.
\end{itemize}
The cardinality of a largest mutual-visibility set, a largest total mutual-visibility set, a largest outer mutual-visibility set, and a largest dual mutual-visibility set will be respectively denoted by $\mu(G)$, $\mut(G)$, $\muo(G)$, and $\mud(G)$. Also, these graph invariants will be respectively called the {\em mutual-visibility number}, the {\em total mutual-visibility number}, the {\em outer mutual-visibility number}, and the {\em dual mutual-visibility number} of $G$. Moreover, for any invariant $\tau(G)$ from the above ones, by a {\em $\tau$-set} we mean any set of vertices of cardinality $\tau(G)$. In addition, for any two invariants $\tau_1(G)$ and $\tau_2(G)$, by $(\tau_1,\tau_2)$-graph we mean any graph $G$ with $\tau_1(G)=\tau_2(G)$.

\subsection{Basic properties}
\label{sec:properties}

If $G$ is a graph, then by definition,
\begin{align}
\mu(G) & \ge \muo(G) \ge \mut(G)\quad {\rm and} \label{eq:1}\\
\mu(G) & \ge \mud(G) \ge \mut(G)\label{eq:2}\,.
\end{align}

In what follows, we first recall some known results about $\mu(G)$, and then use these results to begin giving the reader a first glimpse of the differences among the proposed mutual-visibility variants.

It is easy to observe that $\mu(G)\ge 1$ and $\muo(G)\ge 1$ for each graph $G$ (indeed, any vertex $u\in V(G)$ forms both a mutual-visibility and an outer mutual-visibility set of $G$). Concerning small values of $\mu$, from~\cite[Lemma~4.1]{DiStefano-2022} we know that:
\begin{itemize}
\item $\mu(G)=1$ if and only if $G \cong K_1$;
\item $\mu(G)=2$ if and only if $G \cong P_n$, $n\ge 2$;
\end{itemize}
Moreover, a partial characterization for $\mu(G)=3$ is provided in~\cite{Cicerone-2023}. The following property concerns $(\mu,\mut)$-graphs:

\begin{property} {\rm \cite[Proposition~3.3]{{Cicerone-2022+}}}
Block graphs (and hence trees and complete graphs) and graphs containing a universal vertex are all $(\mu,\mut)$-graphs.
\end{property}

For trees, in~\cite{DiStefano-2022} it is shown that $\mu(T)$ equals the number $L(T)$ of leaves of $T$. Hence, as $T$ is a $(\mu,\mut)$-graph, by~\eqref{eq:1} and~\eqref{eq:2} we get $\tau(T) = L(T)$, where $\tau$ is any of the four mutual-visibility variants. When $T$ reduces to a path $P_n$, we get $\tau(P_n) = 2$. In particular, if $V(P_n) = [n]$, then
\begin{itemize}
\item
 each pair of distinct vertices of $P_n$ forms a $\mu$-set of $P_n$,
\item
 $\{1,n\}$ is the only $\muo$-set of $P_n$,
\item
 $\{1,2\}$, $\{1,n\}$, and $\{n-1,n\}$ are the only $\mud$-sets of $P_n$,
\item
 $\{1,n\}$ is the only $\mut$-set of $P_n$.
\end{itemize}

The latter suggests the following general property: if $G$ is a graph with $n(G)\ge  2$, and $u$ and $v$ form a \emph{diametral pair} -  i.e., $d_G(u,v) = \diam(G)$, then it can be easily verified that $\{u,v\}$ is an outer mutual-visibility set of $G$ (and hence $\muo(G)\ge 2$). As a consequence we get the following:
\begin{itemize}
\item $\muo(G)=1$ if and only if $G \cong K_1$.
\end{itemize}

We have already observed that $\mu(G)\ge 1$ and that $\muo(G)\ge 1$ for each graph $G$. Interestingly, it is known that there exist graphs $G$ such that $\mut(G) = 0$. A characterization of the graphs satisfying this is known as follows.

\begin{theorem}{\em \cite{tian-2023+}}
\label{thm:main-characterization-for-0}
Let $G$ be a graph with $n(G)\ge 2$. Then $\mut(G) = 0$ if and only if each vertex of $G$ is the middle vertex of a convex $P_3$ in $G$.
\end{theorem}

In the same lines of the situation that happens with the total mutual-visibility number, there are also graphs $G$ for which $\mud(G) = 0$. A simple example of this is for instance the cycle $C_7$. We next give some partial results concerning characterizing the graphs achieving this property.

\begin{proposition}
\label{prop:mud-0}
Let $G$ be a graph. If every two adjacent vertices of $G$ are the center of a convex $P_4$, then $\mud(G)=0$.
\end{proposition}

\proof
Let $G$ be a graph in which every two adjacent vertices are the center of a convex $P_4$. We first observe that $\delta(G)\ge 2$ for otherwise the edge between a leaf and its unique adjacent vertex is not the center of a convex $P_4$. Suppose now on the contrary that $\mud(G)\ge 1$ and let $u\in V(G)$ be a vertex that lies in some dual mutual-visibility set $S$. Let $u'$ be a neighbor of $u$. As the edge $uu'$ is the center of a convex $P_4$, there exist vertices $w$ and $w'$ such that $w, u, u', w'$ is a convex $P_4$. Then at least one of $w$ and $u'$ must lie in $S$ for otherwise these two vertices are not $S$-visible. 
Suppose that $w\in S$. Consider a convex $P_4$ such that the edge $wu$ is its center, say $x, w, u, y$ (where it is possible that $y = u'$). If $x\in S$, then $x$ and $u$ are not $S$-visible, and if $x\notin S$, then $x$ and $y$ are not $S$-visible. This contradiction implies that $w\notin S$. Then $u'\in S$, for otherwise $w$ and $u'$ are not $S$-visible. This implies that $w'\notin S$, for otherwise $u$ and $w'$ are not $S$-visible, but then $w$ and $w'$ are not $S$-visible, a final contradiction. 
\qed

Proposition \ref{prop:mud-0} can be used to get a characterization of the graphs $G$ with girth at least $7$ for which $\mud(G)=0$.

\begin{proposition}
\label{prop:mud-girth-7}
Let $G$ be a graph with girth at least $7$. Then $\mud(G)=0$ if and only if $\delta(G) \ge 2$.
\end{proposition}

\proof
First, if $\mud(G)=0$, then we readily observe that $G$ has minimum degree at least two. On the other direction, if $G$ has minimum degree at least $2$, then clearly every two adjacent vertices of $G$ are the center of a convex $P_4$, since $G$ has girth at least $7$. Thus, by Proposition \ref{prop:mud-0} we obtain that $\mud(G)=0$.
\qed

Notice that there are graphs $G$ not satisfying the statement of Proposition \ref{prop:mud-0} such that $\mud(G)=0$. Examples of this are for instance the tori $C_5\cp C_5$, $C_5\cp C_6$ and $C_6\cp C_6$. The proofs of these facts shall be given in Proposition \ref{thm:tori-dual}.

The following statement refers to some useful properties of mutual-visibility parameters (one of them already proved in~\cite{tian-2023+}).

\begin{proposition}
If $X$ is a mutual-visibility set (outer mutual-visibility set, total mutual-visibility set, respectively) of a graph $G$ and $Y\subseteq X$, then $Y$ is also a mutual-visibility set (outer mutual-visibility set, total mutual-visibility set, respectively) of $G$.
\end{proposition}
\begin{proof}
Let $X\subseteq V(G)$, $u\in X$, and $X'=X\setminus \{u\}$. Consider the following incremental hypothesis:
\begin{enumerate}
\item
Assume that $v$ and $w$ are $X$-visible for each $v,w\in X$. Then, trivially, $v$ and $w$ are $X'$-visible for each $v,w\in X'$.
\item
Additionally to the hypothesis in the previous item, assume also that $v$ and $w$ are $X$-visible for each $v\in X$ and for each $w\in \overline{X}$. Hence, $v$ and $w$ are $X'$-visible for each $v\in X'$ and for each $w\in \overline{X'}$ (if $w=u$, the property holds since $X$ is a mutual-visibility set).
\item
Additionally to the hypothesis in the previous two items, assume also that $v$ and $w$ are $X$-visible for each $v,w\in \overline{X}$. Hence, $v$ and $w$ are $X'$-visible for each $v,w\in \overline{X'}$ (if $v=u$ or $w=u$, the property holds since $X$ is an outer mutual-visibility set).
\end{enumerate}
From Item~1, we get the conclusion for the mutual-visibility. From Items~1 and~2, the conclusion is deduced for the outer-visibility. From all the three items above, the total-visibility property is also obtained.
\qed
\end{proof}

Notice that the property of the previous statement does not hold for the dual mutual-visibility. In fact, let $u,v$ be an edge in the cycle $C_6$: it can be observed that $X=\{u,v\}$ is a $\mud$-set of $C_6$, whereas $X'=\{u\}$ is not a dual mutual-visibility set of the same graph (the two adjacent vertices of $u$ are not $X'$-visible).

\smallskip
We end this section by listing the exact values for $\tau(C_n)$ when $\tau$ is any of the four mutual-visibility variants. As for paths, also this special kind of graphs allows us to emphasize the different behavior of the mutual-visibility variants. 

\smallskip
Concerning the original mutual-visibility, from~\cite{DiStefano-2022} we recall that
\begin{equation}\label{eq:mu_Cn}
\mu(C_n) = 3,\ n\ge 3\,.
\end{equation}
%
For the total mutual-visibility we have (cf.~\cite{tian-2023+}):
\begin{equation}\label{eq:mut_Cn}
\mut(C_n)=\left\{\begin{array}{ll}
                      3; & n=3, \\
                      2; & n=4, \\
                      0; & n\ge 5.
                    \end{array}
\right.
\end{equation}
%
Consider next the dual mutual-visibility number of cycles. If $n\ge 7$, then Proposition~\ref{prop:mud-girth-7} gives $\mud(C_n)=0$. We can easily process short cycles so that we have: 
\begin{equation}\label{eq:mud_Cn}
\mud(C_n)=\left\{\begin{array}{ll}
                      3; & n\in \{3,4\}, \\
                      2; & n\in \{5, 6\},\\
                      0; & n\ge 7.
                    \end{array}
\right.
\end{equation}
Finally, for the outer mutual-visibility number we have:
\begin{equation}\label{eq:muo_Cn}
\muo(C_n)=\left\{\begin{array}{ll}
                      3; & n = 3, \\
                      2; & n\ge 4.
                    \end{array}
\right.
\end{equation}
%
The results \eqref{eq:mut_Cn}-\eqref{eq:muo_Cn} will be used in Section~\ref{sec:grids} for determining the values of the corresponding invariants for grids and tori, and in Section~\ref{sec:compare} for comparing the four invariants.

\section{Complexity of the problems}
\label{sec:complexity}

In this section we address computational issues concerning the the three mutual visibility parameters considered in our investigation. 

%

\smallskip
We start remarking that testing whether $X$ is a $\tau$-set, with $\tau\in \{\muo,\mud,\mut\}$, can be performed in polynomial time. In~\cite{DiStefano-2022}, it is shown that checking whether $X\subseteq V(G)$ is a $\mu$-set of a connected graph $G$ can be computed in $O(|X|m(G))$ time. This test is performed as follows: given a vertex $u\in X$, the distance in $G$ between $u$ and any other vertex can be computed by performing a standard BFS in $G$ that starts from $v$; moreover, by using an adapted BFS it is also possible to compute these distances with the constraint that the shortest paths cannot use any element of $X$ as internal vertex. By comparing the distances computed by the two BFS it is possible to verify whether $u$ and any other vertex $v\in X$ are $X$-visible. 
It can be easily observed that the same approach can be used for testing whether $X$ is a $\tau$-set, with $\tau\in \{\muo,\mud,\mut\}$. For instance, testing a $\muo$-set requires two BFS for each vertex $u\in X$ (but now the output of the two BFS also concerns the distances between  $u$ and any other vertex in $V(G)$) - this leads to the same complexity of $O(|X|m(G))$. Instead, for testing a $\muo$-set, two BFS for each vertex $u\in X$ and two BFS for each vertex $u\in \overline{X}$ are required, thus leading to a total complexity of $O(n(G)m(G))$. Of course, the latter complexity also holds for testing a $\mut$-set.


\medskip
In the remainder of the section we show that the decision problems regarding computing the values of the three mutual visibility parameters are computationally difficult. In addition, we remark that the provided reduction can be used to prove the complexity of the standard mutual-visibility problem, which was already proved in \cite{DiStefano-2022}. For such purposes, we define the {\sc $\tau(G)$-mutual-visibility problem}, written in general for any parameter $\tau(G)$ with $\tau(G)\in \{\mu(G),\mud(G),\muo(G),\mut(G)\}$:

\begin{itemize}
\item[-]
{\sc Instance}: A graph $G$, a positive integer $k\leq n(G)$, and a given parameter $\tau(G)\in \{\mu(G),\mud(G),\muo(G),\mut(G)\}$. 
\item[-]
{\sc Question}: Is it satisfied that $\tau(G)\geq k$?
\end{itemize}

\begin{theorem}
For any parameter $\tau(G)\in \{\mu(G),\mud(G),\muo(G),\mut(G)\}$, the {\sc $\tau(G)$-mutual-visibility problem} is NP-complete.
\end{theorem}

\proof
For any $\tau(G)\in \{\mu(G),\mud(G),\muo(G),\mut(G)\}$, the {\sc $\tau(G)$-mutual-visibility problem} is in NP, since we have already observed how it is possible to verify in polynomial time that a given set of cardinality at least $k$ is indeed a set a (dual, outer or total) mutual visibility set.

\begin{figure}[!ht]
\centering
\begin{tikzpicture}[scale=1.5, style=thick]
\def\vr{3pt}
\def\len{1}

\coordinate(G) at (-1,0);
\coordinate(a) at (0,0);
\coordinate(b) at (1,0);
\coordinate(c) at (2,0);
\coordinate(d) at (3,0);
\coordinate(e) at (4,0);
\coordinate(x) at (2,1);
\coordinate(x1) at (0.5,-1);
\coordinate(x2) at (1.5,-1);
\coordinate(x3) at (2.5,-1);
\coordinate(x4) at (3.5,-1);

\coordinate(k1) at (2.4,2);
\coordinate(k2) at (0.7,-2.2);
\coordinate(k3) at (1.7,-2.2);
\coordinate(k4) at (2.7,-2.2);
\coordinate(k5) at (3.7,-2.2);
\coordinate(k6) at (4.2,-1);

\coordinate(l1) at (1.6,1.9);
\coordinate(d1) at (2.4,1.9);
\coordinate(l2) at (0.1,-1.9);
\coordinate(d2) at (0.9,-1.9);
\coordinate(l3) at (1.1,-1.9);
\coordinate(d3) at (1.9,-1.9);
\coordinate(l4) at (2.1,-1.9);
\coordinate(d4) at (2.9,-1.9);
\coordinate(l5) at (3.1,-1.9);
\coordinate(d5) at (3.9,-1.9);

\draw (a)--(b)--(c)--(d)--(e);

\draw (a)--(x);
\draw (b)--(x);
\draw (c)--(x);
\draw (d)--(x);
\draw (e)--(x);

\draw (a)--(x1);
\draw (b)--(x1);
\draw (b)--(x2);
\draw (c)--(x2);
\draw (c)--(x3);
\draw (d)--(x3);
\draw (d)--(x4);
\draw (e)--(x4);

\draw (x)--(l1);
\draw (x)--(d1);
\draw (x1)--(l2);
\draw (x1)--(d2);
\draw (x2)--(l3);
\draw (x2)--(d3);
\draw (x3)--(l4);
\draw (x3)--(d4);
\draw (x4)--(l5);
\draw (x4)--(d5);

\draw(a)[fill=white] circle(\vr);
\draw(b)[fill=white] circle(\vr);
\draw(c)[fill=white] circle(\vr);
\draw(d)[fill=white] circle(\vr);
\draw(e)[fill=white] circle(\vr);
\draw(x)[fill=white] circle(\vr);
\draw(x1)[fill=white] circle(\vr);
\draw(x2)[fill=white] circle(\vr);
\draw(x3)[fill=white] circle(\vr);
\draw(x4)[fill=white] circle(\vr);

\draw (0.5,-2) ellipse (0.4cm and 0.25cm);
\draw (1.5,-2) ellipse (0.4cm and 0.25cm);
\draw (2.5,-2) ellipse (0.4cm and 0.25cm);
\draw (3.5,-2) ellipse (0.4cm and 0.25cm);
\draw (2,2) ellipse (0.4cm and 0.25cm);
\draw (2.1,-1) ellipse (2.1cm and 0.3cm);

\draw[anchor = west] (x)++(0.05,0.05) node {$x$};
\draw[anchor = west] (x1)++(0.05,0.0) node {$v_{e_1}$};
\draw[anchor = west] (x2)++(0.05,0.0) node {$v_{e_2}$};
\draw[anchor = west] (x3)++(0.05,0.0) node {$v_{e_3}$};
\draw[anchor = west] (x4)++(0.05,0.0) node {$v_{e_4}$};
\draw[anchor = west] (G) node {$P_5:$};

\draw[anchor = west] (k1) node {$K_t$};
\draw[anchor = west] (k2)++(0.05,0.0) node {$K_t$};
\draw[anchor = west] (k3)++(0.05,0.0) node {$K_t$};
\draw[anchor = west] (k4)++(0.05,0.0) node {$K_t$};
\draw[anchor = west] (k5)++(0.05,0.0) node {$K_t$};
\draw[anchor = west] (k6) node {$K_m$};

\end{tikzpicture}
\caption{Construction of the graph $G'$ from $G=P_5.$}
\label{fig:G'}
\end{figure}
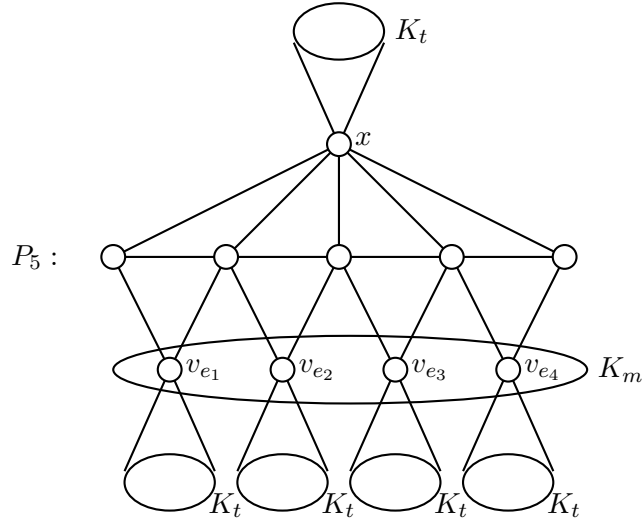

To prove the NP-completeness, we shall make a reduction from the {\sc independent set problem}, which is a classical problem in graph theory known to be NP-complete. We recall that the independence number $\alpha(G)$ of a graph $G$ is the cardinality of a largest edgeless set of vertices of $G$. We consider an arbitrary graph $G$, and will construct a graph $G'$ as follows.
We begin with the graph $G$ of order $n$ with vertex set $V(G) = [n]$. Next for each edge $e=ij$ of $G$, we add an isolated vertex $v_e=v_{ij}$ and the edges $iv_e$ and $jv_e$. Moreover, we add all possible edges between all the vertices $v_e$ with $e\in E(G)$ (namely those vertices form a clique $K_m$). We next add a clique $K_{t+1}$, with $t\ge 3$, and select one of its vertices, denoted by $x$, and join all the vertices of $G$ by an edge to the vertex $x$. Let us say that $V(K_{t+1})=\{x,x_1,\dots,x_t\}$. Finally, for each of the vertices $v_e$ with $e\in E(G)$, we add a clique $K_t$ ($t\ge 3$) with vertex set $V(K_t)=\{e_{y_1},\dots,e_{y_t}\}$ and join all the vertices of such $K_t$ with the corresponding $v_e$. An example of the graph $G'$, for $G = P_5$, is given in Fig.~\ref{fig:G'}.


Let $I$ be an independent set of $G$ of cardinality $\alpha(G)$ and let $S=I \cup \{ x_{i}\,: \, i\in [t]\}\cup \{e_{y_i}: \, i\in [t], e\in E(G)\}$. Hence, clearly $|S|=(m+1)t + \alpha(G)$. We claim that $S$ is total mutual-visibility set of $G'$. We consider some situations for any pair of vertices $u,v\in V(G')$. First note that, if $u,v$ are adjacent, then they are $S$-visible. This includes several possible cases (for instances pairs $(v_e,v_f)$, $(e_{y_i},e_{y_j})$, $(x_i,x_j)$, or $(x,x_j)$ and some other ones). In this regard, from now on we assume that $u$ and $v$ are not adjacent.

\begin{enumerate}
\item If $u=e_{y_i}$ and $v=e'_{y_j}$ for $i,j\in [t]$, $e,e'\in E(G)$, then $u,v$ are $S$-visible since the geodesic $u=e_{y_i}, v_{e}, v_{e'}, e'_{y_j}=v$ does not contain vertices of $S$ (other than $u,v$).
\item If $u=e_{y_i}$ and $v=v_{e'}$ for $i,j\in [t]$, $e,e'\in E(G)$, then the geodesic $u = e_{y_i}, v_{e}, v_{e'}=v$ does not contain vertices of $S$ (other than $u,v$), and so, $u,v$ are $S$-visible.
\item If $u=e_{y_i}$ and $v\in V(G)$ for $i\in [t]$ and $e\in E(G)$, then $u,v$ are $S$-visible either by the geodesic $u=e_{y_i}, v_{e}v$ (when $v\in e$) or by a geodesic $u=e_{y_i}, v_{e}, v_{e'}, v$ (when $v\notin e$ and $v\in e'$).
\item If $u=e_{y_i}$ and $v=x$  for $i\in [t]$ and $e\in E(G)$, then they are $S$-visible through a geodesic $u=e_{y_i}, v_{e}, \ell, x=v$ where $\ell\in e$. Notice that such $\ell$ always exists because $I$ is independent and for each edge $ij\in E(G)$ at least one of $i$ or $j$ is not in $I$.
\item If $u=e_{y_i}$ and $v=x_j$ for $i,j\in [t]$ and $e\in E(G)$, then they are $S$-visible through a geodesic $u=e_{y_i}, v_{e}, \ell, x, x_j=v$ where $\ell\in e$.
\item If $u=v_{e}$ and $v\in V(G)$ for $e\in E(G)$, then they are $S$-visible through a geodesic $u=v_e, v_{e'}, \ell=v$ where $\ell\in e'$.
\item If $u=v_{e}$ and $v=x$ for $e\in E(G)$, then they are $S$-visible through a geodesic $u=v_e, \ell, x=v$ where $\ell\in e$. Again, this $\ell$ always exists because $I$ is independent and for each edge $ij\in E(G)$ at least one of $i$ or $j$ is not in $I$.
\item If $u=v_{e}$ and $v=x_j$ for $j\in [t]$ and $e\in E(G)$, then they are $S$-visible through a geodesic $u=v_{e},\ell, x, x_j=v$ where $\ell\in e$.
\item If $u,v\in V(G)$, then they are $S$-visible through the geodesic $u, x, v$, since they are at distance $2$ (notice that the situation when $u,v$ are adjacent is already mentioned).
\item If $u\in V(G)$ and $v=x_j$ for $j\in [t]$, then they are $S$-visible through the geodesic $u, x, x_j=v$.
\end{enumerate}

These cases cover all possible pairs of vertices of $G'$, which shows that any two vertices of $G'$ are $S$-visible, and so $S$ is a total mutual-visibility set of $G'$.  Therefore, by using \eqref{eq:1} and \eqref{eq:2}, we have that
\begin{equation}\label{eq:eq-3}
\mu(G')\ge \mud(G')\ge \mut(G')\ge |S|=(m+1)t+\alpha(G)
\end{equation}
and
\begin{equation}\label{eq:eq-4}
\mu(G')\ge \muo(G')\ge \mut(G')\ge |S|=(m+1)t+\alpha(G).
\end{equation}

On the other hand, we consider now the parameter $\mu(G')$ and let $S$ be a $\mu$-set of $G'$. We first observe that vertices from all the copies of the complete graphs $K_t$ or from the graph $K_{t+1}$ minus the vertex $x$, used to construct $G'$ are in $S$ because these vertices do not lie on a shortest path between other vertices of $G'$. Furthermore, none of the vertices $x$ and $v_{e}$, with $e \in E(G)$, are in $S$. Indeed, if for instance $x \in S$, then the vertices $x_i\in V(K_{t+1})\setminus\{x\}$, $i\in [t]$ are not visible with any vertex outside $V(K_{t+1})\setminus\{x\}$. The same argument applies for vertices $v_{e}$, with $e \in E(G)$.

Now suppose there exists a mutual-visibility set $S'$ of cardinality at least $(m+1)t + \alpha(G) +1$. Since the vertices $e_{y_i},x_j$, for $i,j\in [t]$ and $e\in E(G)$ are in $S'$ and the vertices $x,v_{e'}$ for $e'\in E(G)$ are not in $S'$, it must happen $|S''|=|S'\cap V(G)|\ge \alpha(G)+1$ ($S''$ is not independent). However, this means that there exist vertices $i,j \in S''$, such that $e=ij \in E(G)$. Hence, every vertex $e_{y_\ell}$ with $\ell\in [t]$ is not $S$-visible with every vertex $x_q$ with $q \in [t]$, since the only two shortest paths between them contain either $i$ or $j$ (the end vertices of $e$). This is a contradiction with $S'$ being a mutual-visibility set. Therefore, $\mu(G') \leq (m+1)t+\alpha(G)$. Again by using  \eqref{eq:1} and \eqref{eq:2} we have that
\begin{equation}\label{eq:eq-5}
(m+1)t+\alpha(G)\ge \mu(G')\ge \mud(G')\ge \mut(G')
\end{equation}
and
\begin{equation}\label{eq:eq-6}
(m+1)t+\alpha(G)\ge \mu(G')\ge \muo(G')\ge \mut(G').
\end{equation}
Now, by using \eqref{eq:eq-3}-\eqref{eq:eq-6} altogether, we deduce that
$$\mu(G')=\muo(G')=\mud(G')=\mut(G')=(m+1)t+\alpha(G),$$
which completes the reduction of the {\sc independent set problem} to the {\sc $\tau(G)$-mutual-visibility problem} where $\tau(G)\in \{\mu(G),\mud(G),\muo(G),\mut(G)\}$.
\qed

\section{Mutual-visibility numbers in grid like structures}
\label{sec:grids}

From the seminal paper on mutual-visibility in graphs we recall the following result.

\begin{theorem} {\rm \cite[Theorem~4.6]{DiStefano-2022}}
\label{thm:mu-grids}
If $n\ge 4$ and $m\ge  4$, then $\mu(P_n\cp P_m) = 2 \min\{n,m\}$.
\end{theorem}

In~\cite[Theorem 4.5]{tian-2023+} it is proved that if $T$ is tree with $n(T)\ge 3$ and $H$ is a graph with $n(H)\ge 2$, then $\mut(T\cp H) = \mut(T) \mut(H)$. As a consequence, by using some induction procedure, we have:

\begin{corollary}
\label{cor:mut-grids}
If $H_k = P_{n_1}\cp \cdots \cp P_{n_k}$, where $k\ge 2$ and  $n_i\ge 3$ for $i\in [k]$, then $\mut(H_k) = 2^k$.
\end{corollary}

In the rest of the section we determine $\muo(G)$, $\mud(G)$, and $\mut(G)$ when $G$ is a grid graph and a torus. Observe that the following property is implicitly assumed in most of the provided proofs: if $L$ is any layer of a grid graph, then since $L$ is isomorphic to a path and is a convex subgraph of the grid, then from $\tau(P_n)=2$ we get that at most two vertices in $L$ can be part of a $\tau$-set of the grid, where $\tau$ is any of the four mutual-visibility variants. The same property holds in the tori, but here the maximum number of vertices per layer depends on the different values of $\tau(C_n)$ as computed in Section~\ref{sec:properties}, cf.~\eqref{eq:mut_Cn}-\eqref{eq:muo_Cn}.

\medskip
We begin with the outer mutual-visibility of grids.

\begin{theorem}
\label{thm:grids-outer}
If $n\ge m\ge 2$, then
$$\muo(P_n\cp P_m) = \begin{cases}
  2; &  m=n=2,\\
  4; & (n,m) \in \{(3,2), (3,3), (4,3), (4,4)\}, \\
  5; & (n,m) \in \{(5,4), (5,5), (6,4)\}, \\
  6; & n=6, m=5, \\
  m+2; & \mbox{otherwise}.
 \end{cases}$$
\end{theorem}

\proof
Let $n\ge m\ge 2$. For simplicity we write $G = P_n\cp P_m$. Let $X$ be a $\muo(G)$-set.

Let $L$ be an arbitrary layer of $G$. Since layers in $G$ are convex subgraphs isomorphic to a path, $|X\cap V(L)|\le 2$. Moreover, if $|X\cap V(L)| = 2$, then the two vertices from the intersection are the end-vertices of $L$.

Having in mind that $m\le n$, we now consider the layers $P_n^{(j)}$, $j\in [m]$. Assume first that $|X\cap V(P_n^{(j)})| = 2$ holds for some $2\le j\le m-1$. Then $|X\cap V(P_m^{(1)})| = 1$ and $|X\cap V(P_m^{(n)})| = 1$. If follows that $|X\cap V(P_n^{(j')})| \le 1$ holds for each $j'\ne j$. We conclude that in this case $|X| \le m+1$. Assume second that $|X\cap V(P_n^{(1)})| = 2$. Hence these two vertices are $(1,1)$ and $(n,1)$. If for every $j\ge 2$ we have $|X\cap V(P_n^{(j)})| \le 1$, then we get the same conclusion. Assume second that for some other $j\ge 2$ we have $|X\cap V(P_n^{(j)})| = 2$. Then by the above argument, $j=m$. By the same argument, $|X\cap V(P_n^{(j')})| \le 1$ holds for each $j' \notin \{1,m\}$. We conclude that $|X| \le m+2$. We have thus proved that $\muo(G) = |X| \le m+2$.

Let first $m=2$. Then for any $n\ge 3$ we clearly have $\muo(G) = 4$. Let next $m = 3$. Then for $n\in \{3,4\}$ we easily get that $\muo(G) = 4$, while for $n \ge 5$ we have $\muo(G) = 5$, where $\{(1,1), (n,1), (3,2), (1,3), (n,3)\}$ is a $\muo(G)$-set. Let $m = 4$. If $n = 4$, then $\muo(G) = 4$, if $n \in \{5, 6\}$, then $\muo(G) = 5$, and if $n\ge 7$, then $\muo(G) \ge 6 = m + 2$. Let next $m = 5$. Then if $n = 5$, we have $\muo(G) = 5$, and if $n = 6$, then $\muo(G) = 6$. For $n\ge 7$ we see that the set $\{(1,1), (n,1), (5,2), (2,3), (4,4), (1,n), (n,5)\}$ is an outer mutual-visibility set. As in general $\muo(G) \le m+2$, this set is a $\muo(G)$-set.
Moreover, $\muo(P_6\cp P_6) = 8$ as $\{(1,1),(1,6),(3,2),(5,3),(2,4),(4,5),(6,1),(6,6) \}$ is a maximal outer mutual-visibility set.

Assume now that $n\ge 7$ and $m\ge 6$. Define a set $X = W \cup A \cup B$, where
\begin{align*}
W & = \{(1,1), (n,1), (1,m), (n,m)\}, \\
A & = \{(2k+1,k+1):\ 1\le k\le \left\lfloor (n-2)/2\right\rfloor\}, \\
B & = \{(2k,\left\lfloor (n-2)/2\right\rfloor+k+1):\ 1\le k\le m - \left\lfloor (n-2)/2\right\rfloor - 2\}.
\end{align*}

\begin{figure}[htb]
\begin{center}
\includegraphics[width=14cm]{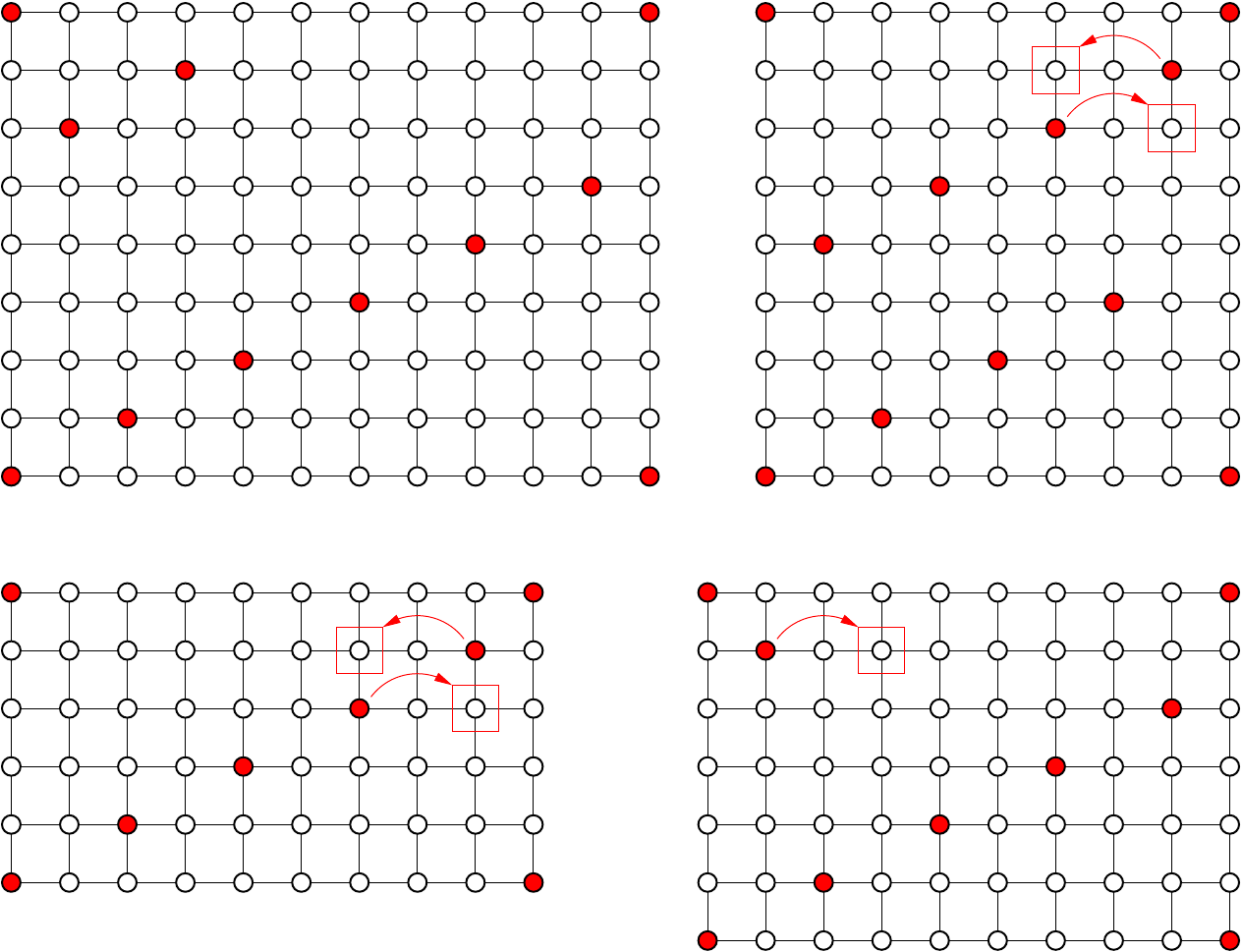}
\caption{Grids $P_{12}\cp P_9$, $P_9\cp P_9$, $P_{10}\cp P_6$, and $P_{10}\cp P_7$ equipped with $\muo$-sets. A red arrow shows how the given $\muo$-set $X$ is transformed into $X'$. }
\label{fig:grids}
\end{center}
\end{figure}

See Fig.~\ref{fig:grids} where the set $X$ is shown for the case $P_{12}\cp P_9$. Notice that the last vertex of $A$ could coincide with the vertex $(n-1,m-1)$. If this happens, we modify $X$ to
$$X' = (X \setminus \{(n-3,m-2), (n-1,m-1)\}) \cup  \{(n-3,m-1), (n-1,m-2)\},$$
see Fig.~\ref{fig:grids}, where the set $X'$ is shown for the case $P_{10}\cp P_6$. Similarly, the last vertex of $B$ could coincide with the vertex $(n-1,m-1)$. If this happens, we also modify $X$ to $X'$, see Fig.~\ref{fig:grids}, where the set $X'$ is shown for the case $P_{9}\cp P_9$. Finally, it can also happen that when $|B|=1$, its vertex can be the vertex $(2,m-1)$. In this case we modify the set $X$ to
$$X'' = (X \setminus \{(2,m-1)\}) \cup  \{(4,m-1)\},$$
see Fig.~\ref{fig:grids}, where the set $X''$ is shown for the case $P_{10}\cp P_7$. Notice that this can be done as $n\ge 7$.

Let $Z$ be either $X$ or $X'$ or $X''$, depending on the case in which we are. Notice first that any two vertices from the set $W = \{(1,1), (n,1), (1,m), (n,m)\}$ are $Z$-visible. Consider now any two vertices $(i,j), (i',j')\in Z$ such that at least one of them is not in $W$. Then $i\ne i'$ and $j\ne j'$. Consider the shortest $(i,j),(i',j')$-path $P$ that lies inside the layers $P_n^{(j)}$ and $P_m^{(i')}$. Then by the construction, $Z\cap V(P) = \{(i,j), (i',j')\}$ and so they are $Z$-visible. We conclude that $Z$ is a mutual-visibility set.

We next show that $Z$ is also an outer mutual-visibility set. Consider vertices $(i,j)\in Z$ and $(i',j')\notin Z$. If $j = j'$ or $i = i'$, then $(i,j)$ and $(i',j')$ are clearly $Z$-visible. Hence assume that $i< i'$ and $j < j'$. The other cases are treated in the same way.

Assume first that $(i,j)\notin W$. If the $(i,j),(i',j')$-geodesic $P$ which simultaneously lies on $P_m^{(i)}$ and $P_n^{(j')}$ has no internal vertex in $Z$, or  the $(i,j),(i',j')$-geodesic $Q$ which simultaneously lies on $P_n^{(j)}$ and $P_m^{(i')}$ has no internal vertex in $Z$, then $(i,j)$ and $(i',j')$ are $Z$-visible. Hence we may assume that both $P$ and $Q$ have an internal vertex each from $Z$. This internal vertex from $P$ cannot lie on $P_m^{(i)}$, hence it is  of the form $(i'',j')$, where $i < i'' < i'$. Similarly, this internal vertex from $Q$ cannot lie on $P_n^{(j)}$, hence it is of the form $(i',j'')$, where $j < j'' < j'$. Moreover, by the construction, at least one of $i'-i'' \ge 2$ and $j'-j'' \ge 2$ holds. We may consider that $i'-i'' \ge 2$ is happening, the other case is done similarly. Then consider the following $(i,j),(i',j')$-geodesic. Take the $(i,j),(i'',j)$-path on $P_n^{(j)}$, proceed with the $(i'',j),(i'',j'-1)$-path on $P_m^{(i'')}$, next take the path $(i'',j'-1)-(i''+1,j'-1)-(i''+1,j')$, and complete the geodesic with the $(i''+1,j'),(i',j')$-path on $P_n^{(j')}$, where we recall that $i''+1 < i'$.

Assume second that $(i,j)\in W$, say $(i,j) = (1,1)$. If $i'< n$ and $j' < m$, then we can argue as above that $(1,1)$ and $(i',j')$ are $Z$-visible. Suppose thus that $(i',j') = (n,j')$, where $j' < m$. In this case consider the following $(1,1),(n,j')$-geodesic. Take the $(1,1),(n-1,1)$-path on $P_n^{(1)}$, proceed with the path $(n-1,1)-(n-1,2)-(n,2)$, and complete the geodesic with the $(n,2),(n,j')$-path on $P_m^{(n)}$.  Hence also in this case $(i,j)=(1,1)$ and $(i',j')=(n,j')$ are $Z$-visible. All the other cases are similar.

We have thus proved that $Z$ is an outer mutual-visibility set. As $|Z| = m+2$, we have proved that $\muo(G) \ge m+2$.
\qed

In the second main result of this section we determine the dual mutual-visibility number of grids.

\begin{theorem}
\label{thm:grids-dual}
If $n\ge 4$ and $m\ge 3$, then $$\mud(P_n\cp P_m) = \left\{\begin{array}{ll}
                                                      3; & n=m=2, \\
                                                      4; & (n,m)=(3,3)\ \mbox{or}\ (n\ge 3\ \mbox{and}\ m=2), \\
                                                      5; & \mbox{otherwise}.
                                                    \end{array}\right.
$$
\end{theorem}

\proof
The case $\mud(P_2\cp P_2) = 3$ is straightforward. Let next assume $n\ge 3$ and $m=2$. Since any two $P_n$-layers can contain at most two vertices from any $\mud$-set, and also because the four vertices $(1,1),(1,2),(n,1),(n,2)$ form a dual mutual-visibility set of $P_n\cp P_2$, we deduce that $\mud(P_n\cp P_2) = 4$. Also, for the case $P_3\cp P_3$ we can check by some simple calculations that also $\mud(P_3\cp P_3) = 4$.

Now on assume $n\ge 4$ and $m\ge 3$ and consider the set $S=\{(1,1), (2,1), (n,m-1), (n,m), (1,m)\}$. It can be readily seen that any two vertices $x,y\in S$ are $S$-visible. Moreover, since $n\ge 4$ and $m\ge 3$, it can be also observed that any two vertices $x',y'\in \overline{S}$ are $S$-visible as well. Thus, $S$ is a dual mutual-visibility set, and so $\mud(P_n\cp P_m) \ge 5$.

To prove the upper bound, let $D$ be a $\mud$-set of $P_n\cp P_m$ (which must have cardinality at least $5$), and consider the following claims that can be easily checked.
\begin{itemize}
  \item Every $P_n$-layer (as well as every $P_m$-layer) in $P_n\cp P_m$ can contain at most two vertices.
  \item If a $P_n$-layer contains two vertices, then at least one of such vertices has the first coordinate from $\{1,n\}$.  Analogous property holds for $P_m$-layers.  
  \item No vertex from the set $\{2,\dots,n-1\}\times \{2,\dots,m-1\}$ could belong to the set $D$. For otherwise, we will find at least two vertices that are not $D$-visible.
  \item No vertex from the sets $\{1,n\}\times \{3,\dots,m-2\}$ and $\{3,\dots,n-2\}\times\{1,m\}$ could belong to the set $D$.
  \item If $(1,1)$, belongs to $D$, then at most one of its neighbors belongs to $D$. For otherwise, since $|D|\ge 5$, the vertex $(1,1)$ would be not $D$-visible with at least two other vertices of $D$.
  \item If $(1,1)$, belongs to $D$ and exactly one of its neighbors, say $(2,1)$, belongs to $D$ too, then no vertex from the sets $\{3,\dots,n\}\times \{1\}$ and $\{1\}\times \{2,\dots,m-1\}$ belongs to $D$.
  \item If $(1,1)$, belongs to $D$ and none of its neighbors belongs to $D$, then $(1,m)$ and $(n,1)$ could belong to $D$ as well, and only such two vertices could satisfy this property in their corresponding layers.
\end{itemize}
As a consequence of the claims above, and up to symmetries, we have the following situations.

\medskip
\noindent
Case 1: $(1,1),(2,1) \in D$. Hence, $(n,1)\notin D$, and this also leads to claim that $(n,2)\notin D$ as well. Since also $\{n\}\times \{3,\dots,m-2\}\cap D=\emptyset$ and $\{1\}\times \{2,\dots,m-1\}=\emptyset$, we deduce that $|D|\le 5$ because $|D\cap ([n]\times\{m\})|\le 2$ and $|D\cap \{(n,m-1)\}|\le 1$.

\medskip
\noindent
Case 2: $(1,1)\in D$ and $(1,2),(2,1) \notin D$. Hence, since $D \cap (\{1\}\times \{2,\dots,m-1\})=\emptyset$ and $D\cap (\{2,\dots,n-1\}\times \{1\})=\emptyset$, then we again obtain that $|D|\le 5$ because $|D\cap ([n]\times\{m\})|\le 2$ and $|D\cap (\{n\}\times [m])|\le 2$.

\medskip
\noindent
Case 3: $(1,1),(1,2),(2,1) \notin D$. If $|D\cap ([n]\times\{1\})|=2$ and $|D\cap (\{1\}\times [m])|= 2$, then $D\cap ([n]\times\{1\})=\{(n-1,1),(n,1)\}$ and $D\cap (\{1\}\times [m])= \{(1,m-1),(1,m)\}$, and it must happen that only the vertex $(n,m)$ could belong to $D$ too. Thus, $|D|\le 5$ in such situation. If without loss of generality $|D\cap ([n]\times\{1\})|=1$ and $|D\cap (\{1\}\times [m])|\in \{1,2\}$, then $D\cap ([n]\times\{1\})=\{(n,1)\}$ and $(1,m)\in D\cap (\{1\}\times [m])$. Thus, we shall obtain that $|D|\le 5$ because $|D\cap ((\{2,\dots,n\}\times \{m\})|\le 1$ and $|D\cap (\{n\}\times \{2,\dots,m\})|\le 1$. Finally, if $|D\cap ([n]\times\{1\})|=0$ and $|D\cap (\{1\}\times [m])|= 0$, then it is straightforward to check that $|D|\le 4$, which is not possible.

\medskip
As a consequence of all the described situations, we deduce that  $\mud(P_n\cp P_m)=|D| \le 5$, which leads to the desired equality and the proof is completed.
\qed

To close the section, we consider the case of the tori. To this end, we need the following proposition.

\begin{proposition}
\label{prop-dual:cover-with-convex}
Let $G$ be a graph. If $V(G) = \bigcup_{i=1}^k V_i$, where $G[V_i]$ is a convex subgraph of $G$ and $\mud(G[V_i]) = 0$ for each $i\in [k]$, then $\mud(G) = 0$.
\end{proposition}
\proof
Suppose on the contrary that $G$ contains a dual mutual-visibility set $X$ with $|X| \ge 1$. Select an arbitrary vertex $x\in X$. Then there exists an $i\in [k]$ such that $x\in V_i$. Hence clearly, $|X\cap G[V_i]| \ge 1$. However, since $G[V_i]$ is convex, we get that $X\cap G[V_i]$ is a dual mutual-visibility set of  $G[V_i]$, a contradiction to the assumption $\mut(G[V_i]) = 0$.
\qed

We recall that, as already observed, if $D$ is a dual mutual-visibility set, a proper subset of $D$ is not necessarily a dual mutual-visibility set, e.g., for $C_n$, $n\in \{5,6\}$, a maximal dual mutual-visibility set is given by two adjacent vertices, whereas a single vertex never represents a dual mutual-visibility set.

\begin{theorem}
\label{thm:tori-dual}
If $n\ge m\ge 3$ then
$$\mud(C_n\cp C_m) = \left\{
 \begin{array}{ll}
  5; & (n,m)\in\{(3,3),(4,3)\}, \\
  8; & (n,m)=(4,4), \\
  2; & (n,m)=(5,3), \\
  4; & (n,m)\in\{(5,4), (6,3),(6,4)\}, \\
  0; & \mbox{otherwise}.
 \end{array}\right.
$$
\end{theorem}

\begin{figure}[ht]
\begin{center}
\includegraphics[width=9.0cm]{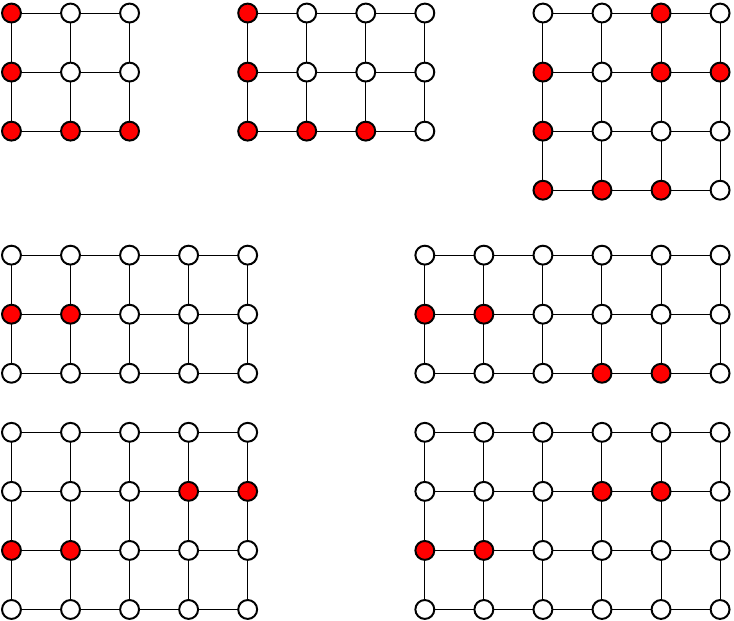}
\caption{Tori, equipped with $\muo$-sets, as defined in the proof of Therorem~\ref{thm:tori-dual}. For sake of simplicity,  each torus is represented as a grid (in each row and column, a wrapping edge is not drawn).}
\label{fig:mud_Cn}
\end{center}
\end{figure}

\proof
Let $V(C_n)=[n]$ and $V(C_m)=[m]$. The different cases in the statement are analyzed as follows (cf.\ Fig.~\ref{fig:mud_Cn}):

\medskip\noindent
Case 1: $n=3, m=3$. It is easy to check that $D=\{(1, 1), (1, 2), (1, 3), (2, 1), (3, 1)\}$ is a dual mutual-visibility set for $C_3\cp C_3$. It is also maximum since a set $D'$ of six or more  vertices cannot be a dual mutual-visibility set. Indeed, if three vertices of $D'$ are on any layer, an adjacent layer has at least two vertices, leading to an induced $C_4$ with the four vertices in $D'$ that are not $D'$-visible. If all the layers have exactly two vertices, then the only possible choice, up to isomorphism, of vertices in $D'$ such that four of them do not induce a $C_4$ implies that the three vertices not in $D'$ are completely surrounded by vertices in $D'$ and then not in dual mutual-visibility with respect to $D'$.

\medskip\noindent
Case 2: $n=4, m=3$. Consider again $D$ and $D'$ as defined in the previous case. Even in this case, $D$ is a dual mutual-visibility while $D'$ is not. Indeed, assume there exists a $C_4$-layer $C_4^{i}$ with three vertices in $D'$. This implies that one among the two $C_4$-layers adjacent to $C_4^{i}$ must have at least two vertices in $D'$. It follows that these vertices in $D'$ induce a $C_4$ cycle or at least two vertices of these two layers not in $D'$ are not $D'$-visible (e.g., if the three vertices in $D'$ of a $C_4$-layer are $(1,1),(2,1),(3,1)$ and the other two vertices are $(3,2),(4,2)$, then $(2,2)$ and $(4,1)$ are not $D'$-visible). Consider now the last case in which all the $C_4$-layers have exactly two vertices in $D'$ and that no four of them induce a $C_4$ cycle. Let $C$ be any $C_4$ cycle of $C_4\cp C_3$ that has two vertices in $D'$. Then these two vertices must be adjacent in $C$ otherwise the remaining vertices in $C$ are not in dual mutual-visibility with respect $D'$. Consequently, the only possible configuration, up to isomorphism, of the vertices in $D'$ is $\{(1,1),(2,1),(2,2),(3,2),(3,3),(4,3)\}$, but then vertices $(2,3)$ and $(3,1)$ are not $D'$-visible. Hence, $D$ is a $\mud$-set.

\medskip\noindent
Case 3: $n=4, m=4$. Let $D^+=D\cup \{(3, 3), (3, 4), (4, 3)\}$. It can be checked that $D^+$ is a dual mutual-visibility set. To show that it is maximal, assume that there exists a dual mutual-visibility set $D'$ of nine or more vertices. Then there is at least one layer with three vertices in $D'$ and an adjacent layer with two adjacent vertices in $D'$. By the discussion in the above case, this is not possible.

\medskip\noindent
Case 4: $(n, m) \in \{(5,3),(5,4)\}$. A $C_5$-layer can have only two adjacent vertices in dual mutual-visibility or none. It can be easily checked that if a $C_5$-layer has two vertices in dual mutual-visibility, then the adjacent layers have no further vertices in dual mutual-visibility. This implies that $\mud(C_5\cp C_3)=2$ and $\mud(C_5\cp C_4)\leq 4$. To see that $\mud(C_5\cp C_4)= 4$ consider the dual mutual-visibility set given by $(1,1),(2,1),(4,3),(5,3)$.

\medskip\noindent
Case 5: $(n, m) \in \{(6,3),(6,4)\}$. As above, a $C_6$-layer can only have two adjacent vertices in dual mutual-visibility or none. Let $D$ be a dual mutual-visibility set for the case $(n,m)=(6,3)$. Without loss of generality, let $(1,2)$ and $(2,2)$ be in $D$. The only vertices that are also in $D$ can be either $(4,1),(5,1)$  or $(4,3),(5,3)$, but not both, showing that $\mud(C_6\cp C_3)= 4$. Similarly for the case $(n,m)=(6,4)$, if $(1,2)$ and $(2,2)$ are in $D$, then only one other $C_6$-layer can have two adjacent vertices in $D$, showing that $\mud(C_6\cp C_4)= 4$.

\medskip\noindent
Case 6: Remaining situations. If $n\geq 7$, then \eqref{eq:mud_Cn} gives $\mud(C_n)=0$ and hence Proposition~\ref{prop-dual:cover-with-convex} implies $\mud(C_n\cp C_m)=0$. The remaining cases are $(n,m)\in\{(5,5),(6,5),(6,6)\}$. Let $D$ a dual mutual-visibility set for any of these cases. Assume, by contradiction, that $D$ is not empty.
We know that $D$ cannot share just one vertex with any $C_n$-layer $C$ (the discussion for a  $C_m$-layer is analogous). Moreover, \eqref{eq:mud_Cn} implies that at most two vertices of $D$ are in $C$. If such two vertices are $(i,j), (i+1,j)$, then the two $C_m$-layers each containing one vertex of $C\cap D$ must both contain an additional vertex, as any $C_5$ or $C_6$ have two vertices in dual mutual-visibility or none. Then one of these pairs of vertices are also in $D$: $\{(i,j+1),(i+1,j+1)\},\{(i,j-1),(i+1,j-1)\},\{(i,j+1),(i+1,j-1)\}$, or $\{(i,j-1),(i+1,j+1)\}$. But each of the first two pairs induces a $C_4$ cycle with $(i,j), (i+1,j)$ vertices in $D$, whereas if $\{(i,j+1),(i+1,j-1)\}$ are in $D$, then $\{(i,j-1),(i+1,j+1)\}$ are not $D$-visible and viceversa.
\qed

\begin{theorem}
\label{thm:tori-total}
If $n\ge m\ge 3$ then
$$\mut(C_n\cp C_m) = \left\{
 \begin{array}{ll}
  3; & (n,m)\in\{(3,3),(4,3)\}, \\
  4; & (n,m)=(4,4), \\
  0; & \mbox{otherwise}.
 \end{array}\right.
$$
\end{theorem}

\proof
It is easy to check that $T=\{(1, 1), (1, 2), (1, 3)\}$ is a total mutual-visibility set for $C_3\cp C_3$. Since any set $T'$ of four vertices implies that there exists an subgraph $C_4$ in $C_3\cp C_3$ with two non adjacent vertices in $T'$ that prevent the visibility of the other two, then $T$ is maximum.

Clearly, the set $T$ is a total mutual-visibility set for $C_4\cp C_3$.  It is also maximum since a set $T'$ of four or more  vertices cannot be a total mutual-visibility set. Indeed, at least two vertices of $T'$, say $u,v$, are on a $C_4$-layer, and to be in total mutual visibility they are adjacent and no further vertex of $T'$ is on the same layer.  Then no vertex of the other two $C_4$-layers can be added to $u,v$ in total mutual visibility. Then $T$ is maximum.

As for $(n,m)=(4,4)$, $\mut(C_n\cp C_m)=4$ and there exists only one $\mut$-set $T$, up to isomorphism. Indeed, if two vertices are on a $C_4$-layer $C$, no vertex of $T$ is on the layers adjacent to $C$. Then at most two other vertices of $T$ can be placed in the fourth $C_4$-layer, for a total of four vertices. A possible $\mut$-set is given by $\{(1,1),(1,2),(3,3),(3,4)\}$.

In all the other cases $m$ is at least $5$ and then, since in this situation $\mut(C_m)=0$, we deduce that $\mut(C_n\cp C_m)=0$.
\qed

The following result provides an upper bound to $\muo(C_n\cp C_m) $ that directly follows from $\muo(C_n)$ as computed in Section~\ref{sec:properties}, cf.~\eqref{eq:muo_Cn}.

\begin{corollary}\label{cor:muo_Cn}
If $n\ge m\ge 3$ then $\muo(C_n\cp C_m) \leq 2m$.
\end{corollary}

\section{Inter-comparison of the mutual-visibility invariants}
\label{sec:compare}

In this section, we compare the four invariants and demonstrate that, roughly speaking, all the theoretical possibilities are feasible.

The results provided in Section~\ref{sec:properties} for cycle graphs indeed demonstrate that the four mutual-visibility numbers can vary pairwise. On the other hand, if $G$ is a $(\mu, \mut)$-graph, then $\mu(G) = \muo(G) = \mud(G) = \mut(G)$. In the above case for cycles, this is demonstrated by the cycle $C_3$. But there are many additional $(\mu, \mut)$-graphs. For instance, in~\cite[Proposition~3.4]{Cicerone-2022+} cographs which are $(\mu, \mut)$-graphs are characterized. In particular, complete split graphs and complete $k$-partite graphs with at least three vertices in each partition set are all $(\mu, \mut)$-graphs.

Next, the results from Section~\ref{sec:grids} show that if $n\ge m\ge 6$, then
$$\mut(P_n\cp P_m) = 4 < \mud(P_n\cp P_m) = 5 < \muo(P_n\cp P_m) = m+2 < \mu(P_n\cp P_m) = 2m\,.$$
This shows that not only can all the four mutual-visibility numbers be pairwise different, but the differences can also be arbitrarily large, except perhaps for $\mut$ and $\mud$. Moreover, from Theorem~\ref{thm:grids-dual} we have $\mud(P_4\cp P_3) = 5 > 4 = \muo(P_4\cp P_3)$, hence $\muo$ and $\mud$ are incomparable.

\begin{figure}[ht!]
\begin{center}
\scalebox{0.8}{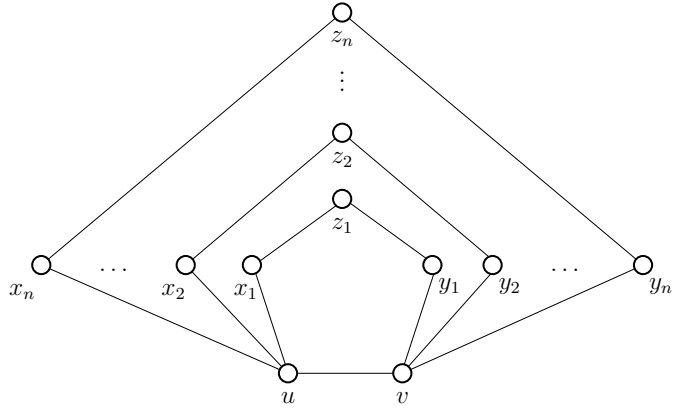}
\caption{The graph $G_n$.}
\label{fig:Gn}
\end{center}
\end{figure}

In the next example we present a family of graphs on which the variety of the four mutual-visibility numbers is again large. The example in addition proves that $\mud$ can be arbitrary larger than $\mut$. Let $G_n$, $n\ge 2$, be the graph obtained from $n$ $5$-cycles having exactly one common edge $uv$. See Fig.~\ref{fig:Gn} from which the definition should be clear and from where the vertex labelling should also be clear.

\begin{proposition}
\label{prop:Gn}
If $n\ge 2$, then $\mu(G_n) = 2n$, $\mud(G_n) = n+1$, $\muo(G_n) = n$, and $\mut(G_n) = 0$.
\end{proposition}

\proof
Let $X$ be a mutual-visibility set of $G_n$. Since $|X\cap V(C)|\le 2$ holds for each $5$-cycle $C$ of $G_n$, we have $\mu(G_n) \le 2n$. On the other hand, the set $\{x_1,\ldots, x_n\} \cup \{y_1,\ldots, y_n\}$ is a mutual-visibility set, hence we can conclude that $\mu(G_n) = 2n$.

Let now $Y$ be an outer mutual-visibility set of $G_n$. Consider the $5$-cycle  $C = ux_1z_1y_1vu$. Note first that no two adjacent vertices of $C$ can belong to $Y$. Suppose now that two non-adjacent vertices of $C$ lie in $Y$. Then we can infer that at least one of the vertices $x_1$, $y_1$, and $z_1$ is $Y$-visible by none of the vertices from $V(G_n)\setminus V(C)$. Hence, in this case we have $|Y| = 2$. In order to have larger outer mutual-visibility sets, each of the cycles can contain at most one vertex. Since $\{z_1, \ldots, z_n\}$ is an outer mutual-visibility set, we conclude that $\muo(G_n) = n$.

The assertion $\mut(G_n) = 0$ follows by an easy verification that $G_n$ fulfills the condition of Theorem~\ref{thm:main-characterization-for-0}.

Next, consider the set $\{x_1, z_1, \ldots, z_n\}$. It can be readily observed that such set is a dual mutual-visibility set of $G_n$, and so $\mud(G_n) \ge n+1\ge 3$ for each $n\ge 2$. Now, let $Z$ be a $\mud$-set of $G_n$. Since $|Z|\ge n+1$, there must be a cycle $C_i = ux_iz_iy_ivu$ with $i\in [n]$ such that $|Z\cap V(C_i)|\ge 2$. Also, by a simple case analysis we infer that $|Z \cap V(C_i)|\le 2$, and so $|Z \cap V(C_i)|=2$. Moreover, in such cycle the vertices of $C_i$ that lie in $Z$ are adjacent, for otherwise the two neighbors (not in $Z$) of any vertex in $Z$ are not $Z$-visible. Clearly, these two adjacent vertices cannot be $u,v$. We may without loss of generality say that $|Z \cap V(C_1)| = 2$.

If $Z \cap V(C_1) = \{x_1,u\}$, then $x_1$ is not $Z$-visible with any other vertex of $G_n$ (except $z_1,y_1$, but such vertices are not in $Z$), and so, $|Z|=2$, which is not possible, since $|Z|\ge 3$. A symmetrical argument also says that $Z \cap V(C_1) \ne \{y_1,v\}$. Thus, by the symmetry of $G_n$, we might assume that $x_1, z_1\in Z$. If there is some $x_i\in Z$ with $i\ne 1$, then such $x_i$ and $z_1$ are not $Z$ visible. Thus, no vertex from $\{x_2,\dots,x_n\}$ is in $Z$. Moreover, for every pair $y_i,z_i$ we have $|Z\cap \{y_i,z_i\}|\le 1$ with $i\in \{2,\dots,n\}$, that is, if $|Z\cap \{y_j,z_j\}|=2$ for some $j$, then $z_j$ and $z_1$  are not $Z$ visible. Now, if there is some $y_j\in Z$ with $j\in \{2,\dots,n\}$, then the two vertices $z_j$ and $v$ (which are not in $Z$) are not $Z$-visible. Thus, no vertex of the set $\{y_2,\dots,y_n\}$ is in $Z$. Consequently, $Z\subset \{x_1,z_1,\dots,z_n\}$, which means $\mud(G_n) \le n+1$, and the proof is completed.
\qed

We close this section by showing that also $\mud$ can be arbitrary larger than $\muo$. To this end, we consider a graph $H_t$, $t\ge 2$, defined as follows. We begin with $t$ copies of the grid graph $P_4\cp P_3$ with vertex sets as previously defined and an extra vertex $x$. Then add, for every copy of the grid $P_4\cp P_3$, the edge between $x$ and the vertex $(2,3)$. Hence $x$ is of degree $t$. 

\begin{proposition}
\label{prop:mud-larger-muo}
If $t\ge 2$, then $\mud(H_t)=5t$ and $\muo(H_t)=4t$.
\end{proposition}

\proof
For each of the copies of $P_4\cp P_3$ in $H_t$, we consider the construction of a $\mud$-set of $P_4\cp P_3$ given in Theorem \ref{thm:grids-dual}. Now, in $H_T$ we consider $S$ as the union of all such sets of vertices from all the copies. It can be readily observed that any two vertices of $S$ are $S$-visible. Also, any two vertices not in $S$ are $S$-visible as well. Thus, $S$ is a dual mutual-visibility set, and so, $\mud(H_t)\ge 5t$. On the other hand, let $D$ is a $\mud$-set of $H_T$. Hence, based on the structure of $H_t$, the distances between vertices from each copy of $P_4\cp P_3$ in $H_t$ are not influenced by the other vertices outside of this copy. Consequently, the restriction of $D$ to any of the copies of $P_4\cp P_3$ in $H_t$ is a dual mutual-visibility set in $P_4\cp P_3$. Therefore, $\mud(H_t)=\sum_{i=1}^{t}|D\cap V(P_4\cp P_3)|\le \mud(P_4\cp P_3)\cdot t=5t$, which leads to the desired equality.

The second formula $\muo(H_t)=4t$ can be obtained by using similar arguments as above, but considering the fact that $\muo(P_4\cp P_3)=4$ as proved in Theorem \ref{thm:grids-outer}.
\qed

Clearly, from Proposition \ref{prop:mud-larger-muo}, we deduce that $\mud$ can be arbitrary larger than $\muo$, which close all the comparisons of the variety of mutual-visibility parameters given in our exposition.

In addition, we might remark that by using analogous techniques to that ones of Proposition \ref{prop:mud-larger-muo}, and other base graph instead of $P_4\cp P_3$, to construct a related $H_t$, some other arbitrarily large (possible) differences between two mutual visibility parameters can be realizable.

\section{Conclusion}\label{sec:conclusion}
%
This work leaves some open problems and suggests some further research directions. Concerning the former, we computed the exact value of $\tau(C_n\cp C_m)$ for each variant $\tau$ except for $\muo$ (for that, Corollary~\ref{cor:muo_Cn} provides just an upper bound). Computing also $\muo(C_n\cp C_m)$ would close the study about tori. Another point is that of providing characterizations for the graphs in which the new invariants have fixed small values (e.g., graphs $G$ for which $\mut(G)=1$, $\mud(G)=0$, or $\mud(G)=1$). For such a task, it could be useful to investigate the notion of \emph{bypass vertices} introduced in~\cite{tian-2023+} as a tool for providing the characterization of Theorem~\ref{thm:main-characterization-for-0}.

We have shown that computing $\tau(G)$ with $\tau(G)\in \{\mu(G),\mud(G),\muo(G),\mut(G)\}$ is an NP-hard problem. It is then worth investigating all the new invariants of our manuscript in special graph classes, with the aim of determining exact formulas for their value or of designing polynomial algorithms for their computation.
Finding structural properties for $\tau$-sets, $\tau\in\{\mu,\muo,\mud,\mut\}$, is also worth to be studied. There is also space for further investigating about the inter-comparison of the mutual-visibility invariants, namely characterizing the graphs achieving equality or being strictly different (smaller or larger) with respect to the values of (some) of the visibility parameters. For instance, as observed in the paper, $\mu(P_n\cp P_m) = 2\cdot \min\{n,m\}$ for $n, m \ge 4$. Also, Theorem~\ref{thm:grids-outer} implies that $\muo(P_n\cp P_m) = \min\{n,m\} + 2$ for almost all $m$ and $n$, which is in the order of $\mu(P_n\cp P_m)/2$. Hence we wonder whether $\mu(G) \le 2\muo(G)$ is true in general. If this is true, then it is sharp by Proposition~\ref{prop:Gn}.

\section*{Acknowledgments}
S. Cicerone and G. Di Stefano were partially supported by the European project  ``Geospatial based Environment for Optimisation Systems Addressing Fire Emergencies'' (GEO-SAFE), contract no. H2020-691161, and by the Italian National Group for Scientific Computation (GNCS-INdAM).
S.\ Klav\v{z}ar was partially supported by the Slovenian Research Agency (ARRS) under the grants P1-0297, J1-2452, and N1-0285.
I.\ G.\ Yero has been partially supported by the Spanish Ministry of Science and Innovation through the grant PID2019-105824GB-I00. Moreover, this investigation was initiated while I.\ G.\ Yero was visiting the University of Ljubljana supported by  ``Ministerio de Educaci\'on, Cultura y Deporte'', Spain, under the ``Jos\'e Castillejo'' program for young researchers (reference number: CAS21/00100).


\end{document}

%% file: Gn.pdf_t
\begin{picture}(0,0)%
\includegraphics{Gn.pdf}%
\end{picture}%
\setlength{\unitlength}{4144sp}%
\begingroup\makeatletter\ifx\SetFigFont\undefined%
\gdef\SetFigFont#1#2#3#4#5{%
  \reset@font\fontsize{#1}{#2pt}%
  \fontfamily{#3}\fontseries{#4}\fontshape{#5}%
  \selectfont}%
\fi\endgroup%
\begin{picture}(4666,3076)(1718,-3505)
\put(3421,-2626){\makebox(0,0)[rb]{\smash{{\SetFigFont{12}{14.4}{\rmdefault}{\mddefault}{\updefault}{\color[rgb]{0,0,0}$x_1$}%
}}}}
\put(3646,-3436){\makebox(0,0)[b]{\smash{{\SetFigFont{12}{14.4}{\rmdefault}{\mddefault}{\updefault}{\color[rgb]{0,0,0}$u$}%
}}}}
\put(4501,-3436){\makebox(0,0)[b]{\smash{{\SetFigFont{12}{14.4}{\rmdefault}{\mddefault}{\updefault}{\color[rgb]{0,0,0}$v$}%
}}}}
\put(2881,-2626){\makebox(0,0)[rb]{\smash{{\SetFigFont{12}{14.4}{\rmdefault}{\mddefault}{\updefault}{\color[rgb]{0,0,0}$x_2$}%
}}}}
\put(1756,-2626){\makebox(0,0)[rb]{\smash{{\SetFigFont{12}{14.4}{\rmdefault}{\mddefault}{\updefault}{\color[rgb]{0,0,0}$x_n$}%
}}}}
\put(2341,-2446){\makebox(0,0)[b]{\smash{{\SetFigFont{12}{14.4}{\rmdefault}{\mddefault}{\updefault}{\color[rgb]{0,0,0}$\ldots$}%
}}}}
\put(5716,-2446){\makebox(0,0)[b]{\smash{{\SetFigFont{12}{14.4}{\rmdefault}{\mddefault}{\updefault}{\color[rgb]{0,0,0}$\ldots$}%
}}}}
\put(4771,-2581){\makebox(0,0)[lb]{\smash{{\SetFigFont{12}{14.4}{\rmdefault}{\mddefault}{\updefault}{\color[rgb]{0,0,0}$y_1$}%
}}}}
\put(5221,-2581){\makebox(0,0)[lb]{\smash{{\SetFigFont{12}{14.4}{\rmdefault}{\mddefault}{\updefault}{\color[rgb]{0,0,0}$y_2$}%
}}}}
\put(6346,-2581){\makebox(0,0)[lb]{\smash{{\SetFigFont{12}{14.4}{\rmdefault}{\mddefault}{\updefault}{\color[rgb]{0,0,0}$y_n$}%
}}}}
\put(4051,-2131){\makebox(0,0)[b]{\smash{{\SetFigFont{12}{14.4}{\rmdefault}{\mddefault}{\updefault}{\color[rgb]{0,0,0}$z_1$}%
}}}}
\put(4051,-1636){\makebox(0,0)[b]{\smash{{\SetFigFont{12}{14.4}{\rmdefault}{\mddefault}{\updefault}{\color[rgb]{0,0,0}$z_2$}%
}}}}
\put(4051,-736){\makebox(0,0)[b]{\smash{{\SetFigFont{12}{14.4}{\rmdefault}{\mddefault}{\updefault}{\color[rgb]{0,0,0}$z_n$}%
}}}}
\put(4051,-1096){\makebox(0,0)[b]{\smash{{\SetFigFont{12}{14.4}{\rmdefault}{\mddefault}{\updefault}{\color[rgb]{0,0,0}$\vdots$}%
}}}}
\end{picture}%